\documentclass[12pt,a4paper]{article}
\usepackage[cp1251]{inputenc}
\usepackage[russian]{babel}
\usepackage[OT1]{fontenc}
\usepackage[left=2.7cm,right=2.7cm,top=3cm,bottom=3 cm]{geometry}
\usepackage[tbtags]{amsmath}
\usepackage{amsfonts,amssymb,mathrsfs,amscd}

\renewcommand{\abstract}[1]{\begin{center}
\begin{minipage}[t]{135mm}
\small
#1
\end{minipage}
\end{center}
}

\makeatletter
\begin{document}

\righthyphenmin=2
\sloppy

\begin{center}
\textbf{{\large A. F. Grishin, Nguyen Van Quynh}}
\bigskip

 \textbf{{\large Entire functions with preassigned zero proximate order }}
\end{center}

\abstract{
It is known that if the proximate order $\rho(r)$ such that $\lim \rho(r)=\rho>0\ (r\to \infty)$, then there exists an entire function $f(z)$ of proximate order $\rho(r)$. In the case where $\rho=0$ the question about the existence of such entire function remained open until now. This question is investigated in the paper.\\
References: 6 units.
\textbf{Comments:} 10 pages, in Russian.\\

\textbf{УДК:} 517.547.22\\
А. Ф. Гришин, Нгуен Ван Куинь, \textbf{Целые функции с наперёд заданным нулевым уточнённым порядком.} Известно, что если уточнённый порядок $\rho(r)$ такой, что $\lim \rho(r)=\rho>0$ $ (r\to \infty)$, то существует целая функция $f(z)$ уточнённого порядка $\rho(r)$. В случае, когда $\rho=0$ вопрос о существовании такой целой функции оставался открытым до настоящего времени. В работе рассматривается вопрос о существовании целой функции уточнённого порядка $\rho(r)$ для случая, когда $\rho=0$.\\
Библиография: 6 названий.\\
\textbf{Ключевые слова:} целая трансцендентная функция, нулевой уточнённый порядок, $\rho$-тригонометрически выпуклая функция.}

Начнём с необходимых определений.

Функция $\rho(r)$ определённая на полуоси $(0,\infty)$ называется уточнённым порядком, если она удовлетворяет условию Липшица на любом сегменте $[a,b]\subset (0,\infty)$ и выполняются условия:\\
\noindent 1) существует предел
$
\rho=\lim\limits_{r\to \infty} \rho(r),
$\\
\noindent 2) $\lim\limits_{r\to \infty} r\rho^{\prime}(r)\ln r=0$ (под $\rho^{\prime}(r)$ следует понимать максимальное по модулю производное число).\\
\noindent В случае, если $\rho=0$, то уточнённый порядок $\rho(r)$ называется нулевым уточнённым порядком.

В работе используется обозначение $V(r)=r^{\rho(r)}$. Сформулируем несколько свойств уточнённого порядка, на которые мы будем ссылаться в дальнейшем.

\medskip
\noindent \textbf{Теорема A. \label{TA}}\textit{ Пусть $\rho(r)$ -- нулевой уточнённый порядок. Тогда для любого $t>0$
$$
\lim\limits_{r\to \infty} \frac{V(rt)}{V(r)}=1,
$$
причём имеет место равномерная сходимость на любом сегменте $[a,b]\subset (0,\infty)$.
}

\medskip
Доказательство можно найти в \cite{book1}, глава 1, $\S$12.

\noindent \rule{80mm}{0.4pt}\\ \smallskip \copyright\hskip2mm
А. Ф. Гришин, Нгуен Ван Куинь, 2014\hfill 

\medskip
\noindent \textbf{Теорема B.\label{TB}}\textit{ 
Пусть $\rho(r)$ -- нулевой уточнённый порядок такой, что $V(r)V(\frac{1}{r})=1.$ Пусть 
$$
\gamma(t)=\sup\limits_{r>0}\frac{V(rt)}{V(r)}.
$$
Тогда $\gamma(t)$ -- непрерывная функция на полуоси $(0,\infty)$, причём функции $\gamma(t)$ и $\gamma(\frac{1}{t})$ имеют нулевой порядок, то есть
$$
\lim\limits_{t\to \infty} \frac{\ln \gamma(t)}{\ln t}=\lim\limits_{t\to \infty} \frac{\ln \gamma(\frac{1}{t})}{\ln t}=0.
$$
}
\medskip
Это теорема 2.5 из \cite{book2}.

\medskip
\noindent \textbf{Теорема C. \label{TC}} \textit{Пусть $\rho(r)$ -- нулевой уточнённый порядок такой, что $V(r)V(\frac{1}{t})=1$. Пусть $K(t)$ -- функция на полуоси $(0,\infty)$ такая, что существует $\delta>0$ такое, что функции $ t^{\pm \delta} K(t)$ принадлежат пространству $L_1(0,\infty)$. Тогда
$$
\lim\limits_{r\to \infty} \frac{1}{rV(r)}\int\limits_{0}^{\infty}K\left(\frac{t}{r}\right)V(t)dt=\int\limits_{0}^{\infty} K(t)dt.
$$
}

Эта теорема легко следует из теорем A и B. Более сильная теорема (теорема 4.1.3) есть в \cite{book3}.

Уточнённые порядки $\rho(r)$ и $\rho_1(r)$ называются эквивалентными, если 
$$
\lim\limits_{r\to \infty}\frac{V(r)}{V_1(r)}=1.
$$

\medskip
\noindent \textbf{Теорема D.\label{TD}} \textit{Пусть $\rho(r)$ -- нулевой уточнённый порядок. Тогда существует эквивалентный ему уточнённый порядок $\rho_1(r)$ такой, что функция $V_1(r)$ допускает голоморфное продолжение в полуплоскость $Re z>0$ и такой, что для любого $n\ge 1$ выполняется равенство
$$
\lim\limits_{r\to \infty} r^n \rho_1^{(n)}(r)\ln r =0.
$$
}

Аналогичные утверждения приведены в \cite{book2}, лемма 2.2 и \cite{book4}, теорема 7. Добавим, что если считать, что функция $V(r)$ удовлетворяет условию $V(r)V(\frac{1}{r})=1$, то в качестве $V_1(z)$ можно взять функцию
$$
V_1(z)=\frac{2z}{\pi}\int\limits_{0}^{\infty}\frac{V(t)}{t^2+z^2}dt,\ Re z>0.
$$

Пусть $f(z)$ -- целая функция. Мы будем пользоваться стандартным обозначением $M(r,f)=\max\limits_{\varphi} |f(re^{i\varphi})|.$

Уточнённый порядок $\rho(r)$ называется уточнённым целой функции $f(z)$, если выполняется условие
$$
\varlimsup\limits_{r\to \infty}\frac{\ln M(r,f)}{V(r)}=\sigma\in (0,\infty).
$$
Если $f(z)$ -- целая функция, то $\upsilon(z)=\ln |f(z)|$ есть субгармоническая функция в комплексной плоскости $\mathbb{C}$.

Как показал Валирон у всякой целой функции конечного порядка $\rho$,
$$
\rho=\varlimsup\limits_{r\to \infty}\frac{\ln \ln M(r,f)}{\ln r},
$$
существует уточнённый порядок $\rho(r)$.

Мы изучаем вопрос:" Какой должен быть уточнённый порядок $\rho(r)$, чтобы он был уточнённым порядком целой трансцендентной функции ?"

Вначале обсудим случай, когда $\lim\limits_{r\to \infty}\rho(r)=\rho >0$ и когда ответ на поставленный выше вопрос известен. 

Нам будет нужно следующее понятие. Дважды дифференцируемая функция $h(\theta)$ на интервале $(\alpha,\beta)$ называется $\rho$-тригонометрически выпуклой, если выполняется неравенство $h^{\prime \prime}(\theta)+\rho^2h(\theta)\ge 0$. В общем случае одно из эквивалентных определений $\rho$-тригонометрической выпуклости таково. Функция $h(\theta)$ называется $\rho$-тригонометрически выпуклой на интервале $(\alpha,\beta)$, если она непрерывна и функция $h^{\prime \prime}(\theta)+\rho^2h(\theta)$, рассматриваемая как обобщённая функция Шварца на интервале $(\alpha,\beta)$, является положительной мерой.

Некоторые свойства $\rho$-тригонометрически выпуклых функций изложены в \cite{book1}, глава 1, $\S$16.

Для целых функций уточнённого порядка $\rho(r)\equiv \rho$ (такие функции называются функциями порядка $\rho$ и нормального типа) Линделёф ввёл понятие индикатора. Это понятие распространяется на произвольные уточнённые порядки и произвольные субгармонические функции. Функция $h(\theta)$ вида
$$
h(\theta)=\varlimsup\limits_{r\to \infty} \frac{\upsilon(re^{i\theta})}{V(r)}
$$
называется индикатором субгармонической функции $\upsilon (z)$ относительно уточнённого порядка $\rho(r)$. Индикатор целой функции $f(z)$ -- это индикатор субгармонической функции $\upsilon(z)=\ln |f(z)|$. 

Известно (\cite{book1}, глава 1, $\S$16), что если $\rho(r)$ -- уточнённый порядок субгармонической функции $\upsilon (z)$, то её индикатор есть ненулевая $\rho$-тригонометрически выпуклая функция. Правда, в \cite{book1} доказательство приведено для случая, когда $\upsilon(z)=\ln |f(z)|$. Однако, оно без изменения распространяется на случай произвольных субгармонических функций.

В пятидесятые годы двадцатого века был открытым вопрос (см. \cite{book1}, глава 2, $\S$1, теорема 3) о существовании целой функции с заданным индикатором относительно заданного уточнённого порядка $\rho(r)$. В дальнейшем этот вопрос получил положительное решение. Теперь можно предложить такое решение. Пусть $h(\theta)$ -- произвольная $2\pi$-периодическая $\rho$-тригонометрически выпуклая функция. Тогда функция $\upsilon(re^{i\theta})=r^{\rho}h(\theta)$ есть субгармоническая в $\mathbb{C}$ функция с индикатором $h(\theta)$ относительно уточнённого порядка $\rho(r)\equiv \rho$. Из этого и теоремы 4 из \cite{book5} следует, что если $\rho(r)$ -- произвольный уточнённый порядок такой, что $\lim\limits_{r\to \infty}\rho(r)=\rho$, то существует субгармоническая функция с индикатором $h(\theta)$ относительно уточнённого порядка $\rho(r)$. Далее применяется теорема о приближении произвольной субгармонической функции $\upsilon(z)$ субгармоническими функциями вида $\ln |f(z)|$, где $f(z)$ -- целая функция. Теоремы такого типа доказывались различными авторами. Основопологающий результат в этом направлении принадлежит Юлмухаметову \cite{book6} .

\medskip
\noindent \textbf{Теорема E.\label{TE}} \textit{
Пусть $\upsilon$ -- субгармоническая функция в $\mathbb{C}$ конечного порядка $\rho>0$. Тогда существует целая функция $f(z)$ порядка $\rho$ такая, что 
$$
\upsilon(z)-\ln |f(z)|=O(\ln |z|),
$$
когда $z$ стремится к бесконечности вне множества кругов $|z-z_j|< r_j$ таких, что $\sum\limits_{r_j\ge r}r_j= o(r^{\rho-\alpha})$ для любого наперёд заданного $\alpha\ge \rho$.
}
\medskip

Из сказанного следует, что для любой $\rho$-тригонометрически выпуклой функции $h(\theta)$ и любого уточнённого порядка $\rho(r)$ такого, что $\lim\limits_{r\to \infty} \rho(r)=\rho>0$ существует целая функция $f(z)$ с индикатором $h(\theta)$ относительно уточнённого порядка $\rho(r)$.

Так как 
$$
\varlimsup\limits_{r\to \infty} \frac{\ln M(r,f)}{V(r)}=\max h(\theta),
$$
то тем самым доказано, что для любого уточнённого порядка $\rho(r)$ такого, что $\lim\limits_{r\to \infty} \rho(r)=\rho>0$ существует целая функция уточнённого порядка $\rho(r)$.

Для случая, когда $\rho=0$ вопрос о существовании такой целой функции оставался открытым до настоящего времени. Мы доказываем следующую теорему.

\medskip
\noindent \textbf{Теорема.} \textit{
Пусть $\rho(r)$ -- заданный нулевой уточнённый порядок. Для того, чтобы существовала целая трансцендентная функция $f(z)$, для которой уточнённый порядок $\rho(r)$ был бы уточнённым порядком этой функции, необходимо и достаточно, чтобы выполнялось равенство
\begin{equation}
\label{Kuin-1}
\lim\limits_{r\to \infty} \frac{\ln r}{V(r)}=0.
\end{equation}
}

\textit{Доказательство.} Н е о б х о д и м о с т ь. Эта часть теоремы очевидна. Действительно, пусть $f(z)$ -- целая трансцендентная функция уточнённого порядка $\rho(r)$. Тогда выполняются следующие соотношения:
\begin{equation*}
\begin{split}
\varlimsup\limits_{r\to \infty}&\frac{\ln M(r,f)}{V(r)} =\sigma \in (0,\infty),\hspace*{1cm} \lim\limits_{r\to \infty}\frac{\ln r}{\ln M(r,f)}=0,\\
& \varlimsup\limits_{r\to \infty} \frac{\ln r}{V(r)}\le \varlimsup\limits_{r\to \infty}\frac{\ln r}{\ln M(r,f)} \varlimsup\limits_{r\to \infty}\frac{\ln M(r,f)}{V(r)}=0.
\end{split}
\end{equation*}
Тем самым равенство (\ref{Kuin-1}) доказано.

Д о с т а т о ч н о с т ь. Для того, чтобы доказать теорему, достаточно её доказать для какого либо уточнённого порядка $\rho_1(r)$ эквивалентного $\rho(r)$. Тогда, как следует из теоремы D, можно дополнительно считать, что $V(r)$ -- дважды непрерывно дифференцируемая функция на полуоси $(0,\infty)$. Из сказанного также следует, что теорему достаточно доказывать для какого-нибудь уточнённого порядка $\rho_1(r)$, который совпадает с $\rho(r)$ на некоторой полуоси $[a,\infty)$. Поэтому имея ввиду (\ref{Kuin-1}), можно дополнительно считать, что на полуоси $[1,\infty)$ выполняется неравенство $V(r)\ge \ln r+1$, а на сегменте $[1,e]$ выполняется равенство $V(r)=\ln r+1$.

Поэтому дальнейшее доказательство теоремы мы будем вести, считая что выполняются дополнительные условия:\\
\noindent 1) $V(r)$ -- дважды непрерывно дифференцируемая функция на полуоси $(0,\infty)$,\\
\noindent 2) при $r\ge 1$ выполняется неравенство $V(r)\ge \ln r+1$,\\
\noindent 3) на сегменте $[1,e]$ выполняется равенство $V(r)=\ln r+1$.

Так как доказательство достаточно длинное, мы разобьём его на несколько этапов.

\medskip
\noindent \textbf{1.} Обозначим $\varphi(x)=V(e^x)$. Функция $\varphi(x)$ обладает свойствами:\\
\noindent 1) $\varphi(x)$ -- дважды непрерывно дифференцируемая функция на полуоси $[0,\infty)$,\\
\noindent 2) при $x\ge 0$ выполняется неравенство $\varphi(x)\ge x+1$,\\
\noindent 3) на сегменте $[0,1]$ выполняется равенство $\varphi(x)=x+1$,\\
\noindent 4) $\lim\limits_{x\to \infty} \frac{\varphi(x)}{x}=\infty$,\\
\noindent 5) $\lim\limits_{x\to \infty}\frac{\varphi^{\prime}(x)}{\varphi(x)}=0$.

Заметим, что свойство 4) следует из равенства (\ref{Kuin-1}), а свойство 5) следует из равенства
$$
\frac{\varphi^{\prime}(x)}{\varphi(x)}=\frac{e^xV^{\prime}(e^x)}{V(e^x)}=\rho(e^x)+xe^x\rho^{\prime}(e^x).
$$

Пусть $k$ -- произвольное вещественное число. Выполняется равенство
$$
\lim\limits_{x\to \infty} \left(\varphi(x)-kx\right)=\infty.
$$
Поэтому функция $\varphi(x)-kx$ на полуоси $[0,\infty)$ имеет минимальное значение, которое мы обозначим через $b(k)$. Заметим, что $b(k)$ -- убывающая непрерывная функция и что $b(k)=1$ при $k\le 1$.

\medskip
\noindent \textbf{2.} Множество решений уравнения $\varphi(x)=kx+b(k)$ есть непустой компакт $F_k$. Обозначим
\begin{equation}
\label{Kuin-2}
\tilde{x}_k=\min F_k,\qquad x_k=\max F_k.
\end{equation}
Прямая $\ell(k):\ y=kx+b(k)$ есть опорная прямая к кривой $\mathscr{L}:\ y=\varphi(x)$ и точки $(\tilde{x}_k,\varphi(\tilde{x}_k)),\ (x_k,\varphi(x_k))$ являются точками опоры для прямой $\ell(k)$ и кривой $\mathscr{L}$.

Тем самым доказано, что у кривой $\mathscr{L}$ есть опорная прямая $\ell(k)$ с угловым коэффициентом $k$ и что кривая $\mathscr{L}$ располагается над прямой $\ell(k)$.

Из свойства 4) функции $\varphi(x)$ следует, что кривая $\mathscr{L}$ не лежит ни в какой полуплоскости $y\le kx+b$. Поэтому множество опорных прямых к кривой $\mathscr{L}$, которые имеют угловой коэффициент $k$, содержит единственный элемент $\ell(k)$.

Мы будем рассматривать величины $\tilde{x}_k,x_k$ как функции переменной $k$. Изучим свойства этих функций. Рассмотрим прямые $\ell(k)$ и $\ell(k_1)$ где $k<k_1$. Точка пересечения этих прямых есть точка $M(\xi,k\xi+b(k))$, где $\xi=\frac{b(k)-b(k_1)}{k_1-k}$.

Пусть $h_1$ -- открытый луч, лежащий на прямой $\ell(k):\ y=kx+b(k)$, который есть образ полуоси $(\xi,\infty)$. Поскольку кривая $\mathscr{L}$ лежит над прямой $\ell(k_1)$, а луч $h_1$ лежит строго под прямой $\ell(k_1)$, то на луче $h_1$ нет точек кривой $\mathscr{L}$. Точка $(x_k,kx_k+b(k))$ лежит на кривой $\mathscr{L}$ и прямой $\ell(k)$ и не лежит на луче $h_1$. Значит $x_k\le \xi$. Аналогично доказывается, что $\xi\le \tilde{x}_{k_1}$. Тем самым доказано неравенство $x_k\le \tilde{x}_{k_1}$.

Так как $\mathscr{L}$ -- гладкая кривая, то опорная прямая к этой кривой является касательной к этой кривой в точке опоры. Поэтому $\varphi^{\prime}(x_k)=k,\ \varphi^{\prime}(\tilde{x}_{k_1})=k_1>k,\ x_k<\tilde{x}_{k_1}$. Из этого неравенства и очевидного неравенства $\tilde{x}_k\le x_k$ следует, что функции $\tilde{x}_k,\ x_k$ являются строго возрастающими. Докажем ещё, что
\begin{equation}
\label{Kuin-3}
\lim\limits_{k\to \infty}x_k=\infty.
\end{equation}
Если это не так, то существует число $a>0$ такое, что $\lim x_k=a\ (k\to \infty)$. Это равенство и равенство $\varphi^{\prime}(x_k)=k$, противоречат тому, что функция $\varphi$ является дважды непрерывно дифференцируемой функцией. Тем самым равенство (\ref{Kuin-3}) доказано.

\medskip
\noindent \textbf{3.} Пусть $g_1$ -- открытый отрезок на прямой $\ell(k)$, являющийся образом интервала $(0,\tilde{x}_k)$ при отображении $y=kx+b(k),\ g_2$ -- открытый луч на прямой $\ell(k)$, являющийся образом луча $(x_k,\infty)$ при том же отображении.

Докажем, что через каждую точку отрезка $g_1$ проходит опорная прямая к кривой $\mathscr{L}$ с некоторым угловым коэффициентом $k_1<k$, а через каждую точку луча $g_2$ проходит опорная прямая к кривой $\mathscr{L}$ с некоторым угловым коэффициентом $k_2>k$.

Пусть $M$ -- произвольная точка луча $g_2$. Её координаты имеют вид $(t,kt+b(k))$, где $t$ -- некоторые число на полуоси $(x_k,\infty)$. Покажем, что существует $\delta=\delta(t)>0$ такое, что при любом $\lambda\in [k,k+\delta)$ луч $A_{\lambda}$ прямой $y-kt-b(k)=\lambda(x-t)$, являющийся образом полуоси $[t,\infty)$ не пересекается с кривой $\mathscr{L}$. Если это не так, то существует последовательность $\xi_n\ge t$ такая , что выполняется равенство
\begin{equation}
\label{Kuin-4}
\varphi(\xi_n)=kt+b(k)+\left(k+\frac{1}{n}\right)(\xi_n-t). 
\end{equation}
Из того равенства и свойства 4) функции $\varphi(x)$ следует, что последовательность $\xi_n$ ограничена. Поэтому у последовательности $\xi_n$ есть сходящаяся подпоследовательность $\xi_{n_m}$. Если $\xi=\lim \xi_{n_m}\ (m\to \infty)$, то $\xi\ge t>x_k$. Беря в равенстве (\ref{Kuin-4}) $n=n_m$ и переходя к пределу при $m\to \infty$, получим равенство $\varphi(\xi)=k\xi+b(k)$. Это равенство противоречит определению $x_k$. 

Тем самым доказано, что существуют числа $\delta>0$ такие, что для любого $\lambda\in [k,k+\delta)$ луч $A_{\lambda}$ не пересекается с кривой $\mathscr{L}$. Пусть $\delta_2$ -- точная верхняя грань таких $\delta$. Очевидно, что $\delta_2<\infty$. Пусть $k_2=k+\delta_2$. Легко проверить, что прямая $y-kt-b(k)=k_2(x-t)$, проходящая через точку $M$ является опорной прямой к кривой $\mathscr{L}$. 

Тем самым часть сформулированного выше утверждения, относящаяся к лучу $g_2$ доказана.

Аналогично доказывается другая часть утверждения относящаяся к отрезку $g_1$. На этом мы заканчиваем доказательство сформулированного нами утверждения.

\medskip
\noindent \textbf{4.} Определим на полуоси $[0,\infty)$ функцию
$$
\varphi_1(x)=\sup\limits_{k}(kx+b(k)).
$$
Для любых $x\ge 0$ и $k\in (-\infty,\infty)$ выполняется неравенство $kx+b(k)\le \varphi(x)$. Из этого следует, что
\begin{equation}
\label{Kuin-5}
\varphi_1(x)\le \varphi(x),\ x\in [0,\infty).
\end{equation}

Так как $
\sup\limits_{k}(kx+b(k))\ge x+b(1)=x+1,$ то выполняется неравенство $\varphi_1(x)\ge x+1$. Так как на сегменте $[0,1]$ выполняется равенство $\varphi(x)=x+1$, то на том же сегменте выполняется равенство $\varphi_1(x)=x+1$. Так как при $k\le 1$ выполняются соотношения $kx+b(k)=kx+1\le x+1$, то справедливо равенство
\begin{equation}
\label{Kuin-6}
\varphi_1(x)=\sup\limits_{k\ge 1}(kx+b(x)).
\end{equation}

Так как при $k\ge 1$ функция $kx+b(k)$ возрастает по переменной $x$, то $\varphi_1(x)$ -- возрастающая функция на полуоси $[0,\infty)$. Так как $\varphi_1(x)$ есть верхняя огибающая семейства линейных функций, то $\varphi_1(x)$ -- выпуклая функция на полуоси $[0,\infty)$.

\medskip
\noindent \textbf{5.} Прямая $\ell(k)$ является опорной прямой для кривой $\mathscr{L}$. Пусть $(\tau,\varphi(\tau))$ -- произвольная точка опоры для прямой $\ell(k)$. Тогда выполняется равенство $k\tau+b(k)=\varphi(\tau)$. Так как $\varphi_1(x)\le \varphi(x)$, то $\varphi_1(\tau)\le \varphi(\tau)=k\tau+b(k)$. С другой стороны из определения $\varphi_1$ следует неравенство $\varphi_1(\tau)\ge k\tau+b(k)=\varphi(\tau)$. Таким образом получаем равенство $\varphi_1(\tau)=\varphi(\tau)$.

Опорная прямая к гладкой кривой является касательной к этой кривой в точке опоры. Поэтому $\varphi^{\prime}(\tau)=k$.

Имеем систему условий:\\
\noindent 1) $\varphi_1(x)\le \varphi(x),$\\
\noindent 2) $\varphi_1(\tau)=\varphi(\tau),$\\
\noindent 3) $\varphi$ -- дифференцируемая функция,\\
\noindent 4) $\varphi_1$ -- выпуклая функция.\\
Из этих условий следует, что функция $\varphi_1$ является дифференцируемой в точке $\tau$ и что выполняется равенство $\varphi^{\prime}_1(\tau)=\varphi^{\prime}(\tau)$.

В частности, функция $\varphi_1$ дифференцируема в точках $\tilde{x}_k$ и $x_k$ и выполняются равенства $\varphi_1^{\prime}(\tilde{x}_k)=\varphi_1^{\prime}(x_k)=k$. Если выполняется неравенство $\tilde{x}_k<x_k$, то на сегменте $[\tilde{x}_k,x_k]$ функция $\varphi_1$ является линейной. Таким образом функция $\varphi_1$ является дифференцируемой на интервале $(\tilde{x}_k,x_k)$. Дифференцируемость $\varphi_1$ в точках $\tilde{x}_k$ и $x_k$ была доказана ранее. Таким образом функция $\varphi_1$ дифференцируема в каждой точке сегмента $[\tilde{x}_k,x_k]$. Здесь уже на необходимости предполагать, что этот сегмент невырожденный.

\medskip
\noindent \textbf{6.} Теперь докажем равенство 
\begin{equation}
\label{Kuin-7}
\lim\limits_{x\to \infty}\frac{\varphi_1(x)}{x}=\infty.
\end{equation}
Из равенства $\varphi_1(x_k)=\varphi(x_k)$ и свойства 4) функции $\varphi(x)$ (см. пункт 1.) следует, что функция $\frac{\varphi_1(x)}{x}$ является неограниченной на полуоси $[1,\infty)$. Поскольку $\varphi_1$ -- выпуклая функция, то $(\varphi_1)_+^{\prime}(x)$ -- возрастающая функция. Поэтому существует $\lim (\varphi_1)_+^{\prime}(x)\ (x\to \infty)$. Если этот предел конечный, то это противоречит неограниченности функции $\frac{\varphi_1(x)}{x}$ на полуоси $[1,\infty)$. Поэтому $\lim \left(\varphi_1\right)_+^{\prime}(x)=\infty \ (x\to \infty)$. Из этого следует равенство (\ref{Kuin-7}).

\medskip
\noindent \textbf{7.} Из равенства (\ref{Kuin-7}) следует, что кривая $\mathscr{L}_1:\ y=\varphi_1(x),\ x\in [0,\infty)$ не лежит ни в какой полуплоскости $y\le kx+b$. Поэтому при любом вещественном $k$ у кривой $\mathscr{L}_1$ имеется не более чем одна опорная прямая с угловым коэффициентом $k$. Прямая $\ell(k)$ есть единственная опорная прямая к кривой $\mathscr{L}$, которая имеет угловой коэффициент $k$. Эта прямая также является опорной к кривой $\mathscr{L}_1$, причём точки $(\tilde{x}_k,\varphi_1(\tilde{x}_k))$ и $(x_k,\varphi_1(x_k))$ являются точками опоры для этой прямой. Эта прямая является единственной опорной прямой к кривой $\mathscr{L}_1$, которая имеет угловой коэффициент $k$. Следовательно, у кривых $\mathscr{L}$ и $\mathscr{L}_1$ одно и тоже множество опорных прямых.

\medskip
\noindent \textbf{8.} Теперь докажем, что функция $\varphi_1(x)$ является дифференцируемой функцией на полуоси $(0,\infty)$. Пусть $\tau$ -- произвольное строго положительное число. Нам нужно доказать существование $\varphi_1^{\prime}(\tau)$. Поскольку $\varphi(x)=x+1$ на сегменте $[0,1]$, то можно считать, что $\tau\ge 1$. Если $\varphi_1(\tau)=\varphi(\tau)$, то как следует из рассуждений, приведённых на этапе 5 доказательства, в этом случае функция $\varphi_1$ будет дифференцируемой в точке $\tau$.

Предположим теперь, что выполняется неравенство $\varphi_1(\tau)<\varphi(\tau)$. Проведём через точку $(\tau,\varphi_1(\tau))$ прямую $\ell_1$ опорную к кривой $\mathscr{L}_1$. Так как функция $\varphi_1$ выпукла, то такая прямая существует. Из рассуждений, приведённых на этапе 7 следует, что прямая $\ell_1$ совпадает с одной из прямой $\ell(k)$. 

Докажем, что выполняется неравенство $\tilde{x}_k\le \tau \le x_k$. Предположим, что выполняется неравенство $x_k< \tau$. Пусть $\tau_1$ такое число, что $x_k<\tau_1<\tau$. Из рассуждений приведённых на этапе 3 следует, что через точку $(\tau_1,k\tau_1+b(k))$, лежащую на прямой $\ell(k)$ можно провести прямую $\ell_2$ опорную к кривой $\mathscr{L}$, угловой коэффициент которой $k_1$ строго больше чем $k$. Прямая $\ell_2$ будет опорной прямой также к кривой $\mathscr{L}_1$. Поэтому кривая $\mathscr{L}_1$ лежит над прямой $\ell_2$. С другой стороны точка $(\tau,k\tau+b(k))$ лежит на кривой $\mathscr{L}_1$ и строго под прямой $\ell_2$. Полученное противоречие доказывает неравенство $\tau\le x_k$. Аналогично доказывается неравенство $\tilde{x}_k\le \tau$. Из рассуждений, приведённых на этапе 4 следует дифференцируемость функции $\varphi_1$ в точке $\tau$. Дифференцируемость функции $\varphi_1$ доказана .

\medskip
\noindent \textbf{9.} Мы заканчиваем изучение свойств функции $\varphi_1$ доказательством равенства
\begin{equation}
\label{Kuin-8}
\lim\limits_{x\to \infty}\frac{\varphi_1^{\prime}(x)}{\varphi_1(x)}=0.
\end{equation}

Обозначим $F=\{x\ge 0:\ \varphi_1(x)=\varphi(x)\}$. Это неограниченное множество.

Пусть $\varepsilon$ -- произвольное строго положительное число. Существует $\xi \in F$ такое, что при $x\ge \xi$ будет выполняться неравенство $\frac{\varphi^{\prime}(x)}{\varphi(x)}<\varepsilon$. Пусть $x>\xi$. Если $x\in F$, то $\frac{\varphi^{\prime}_1(x)}{\varphi_1(x)}=\frac{\varphi^{\prime}(x)}{\varphi(x)}<\varepsilon$. Пусть теперь $x\notin F$. Из рассуждений, приведённых на предыдущем этапе следует, что существует $k$ такое, что выполняется неравенство $\tilde{x}_k\le x \le x_k$. Из этого неравенства и соотношения $x\notin F$ следует неравенство $\tilde{x}_k < x < x_k$. Предположим, что выполняется неравенство $\tilde{x}_k \ge \xi$. Используя то, что $\varphi_1$ -- линейная возрастающая функция на сегменте $[\tilde{x}_k,x]$ получим
$$
\frac{\varphi^{\prime}_1(x)}{\varphi_1(x)}\le \frac{\varphi_1^{\prime}(\tilde{x}_k)}{\varphi_1(\tilde{x}_k)}= \frac{\varphi^{\prime}(\tilde{x}_k)}{\varphi(\tilde{x}_k)}<\varepsilon.
$$
Если же выполняется неравенство $\tilde{x}_k< \xi$, то аналогичное рассуждение нужно применить к сегменту $[\xi,x]$. Тем самым равенство (\ref{Kuin-8}) доказано.

\medskip
Пусть числовая функция $y=\psi(x)$ имеет своей областью определения множество $E$. Надграфик $epi \ \psi$ функции $\psi$ определяется следующим образом
$$
epi \ \psi = \left\{(x,y)\in \mathbb{R}^2:\ x\in E,\ y\ge \psi(x)\right\}.
$$
Можно проверить, что построенная функция $\varphi_1(x)$ является наибольшей выпуклой минорантой функции $\varphi(x)$ и что множество $epi\ \varphi_1$ есть выпуклая оболочка множества $epi\ \varphi$.

\medskip
\noindent \textbf{10.} Функция $\varphi_1(x)$ выпукла на полуоси $[0,\infty)$ , причём $\varphi_1(0)=1$, $\left(\varphi_1\right)^{\prime}_+(0)=1$. Поэтому функцию $\varphi_1(x)$ можно продолжить как выпуклую на вещественную ось $(-\infty,\infty)$, причём так чтобы выполнялись соотношения $$
\lim\limits_{x\to -\infty}\varphi_1(x)=m>0,\ \lim\limits_{x\to -\infty}\varphi_1^{\prime}(x)=0.$$

Обозначим $r^{\rho_1(r)}=V_1(r)=\varphi_1(\ln r)$. Заметим, что функция $\rho_1(r)$ -- это нулевой уточнённый порядок, функция $V_1(r)$ является логорифмически выпуклой, причём выполняются соотношения
$$
\lim\limits_{r\to +0}V_1(r)=m>0,\ \lim\limits_{r\to +0}rV_1^{\prime}(r)=0,\ \lim\limits_{r\to \infty} rV_1^{\prime}(r)=\infty.
$$
Последнее равенство есть следствие равенства (\ref{Kuin-7}) и того, что $rV_1^{\prime}(r)$ есть возрастающая функция.

Далее обозначим
$$
n(r)=\left[rV_1^{\prime}(r)\right],\ N(r)=\int\limits_{0}^{r} \frac{n(t)}{t}dt.
$$
функция $n(r)$ является возрастающей ступенчатой неограниченной функцией с единичными скачками. Выполняются соотношения $n(r)\thicksim rV_1^{\prime}(r)\ (r\to \infty)$, 
$$
\lim\limits_{r\to \infty}\frac{N(r)}{V_1(r)}=1,\ \frac{rN_+^{\prime}(r)}{N(r)}=\frac{n(r)}{N(r)}\to 0\ (r\to \infty).
$$

Из последнего равенства следует, что справедливо равенство $N(r)=r^{\rho_2(r)}$, где $\rho_2(r)$ -- нулевой уточнённый порядок.

Пусть $a_n$ -- точки роста функции $n(r)$,
$$
f(z)=\prod\limits_{n=1}^{\infty}\left(1-\frac{z}{a_n}\right).
$$
Имеем
\begin{equation*}
\begin{split}
\ln M(r,f)&=\sum\limits_{n=1}^{\infty}\ln \left( 1+\frac{r}{a_n}\right) = \int\limits_{0}^{\infty} \ln \left(1+\frac{r}{t}\right) dn(t)\\
&=\int\limits_{0}^{\infty} \frac{rn(t)}{t(t+r)}dt = \int\limits_{0}^{\infty} \frac{rN(t)}{(t+r)^2}dt.
\end{split}
\end{equation*} 
Теперь из теоремы C следует, что 
$$
\lim\limits_{r\to \infty}\frac{\ln M(r,f)}{N(r)}=\int\limits_{0}^{\infty}\frac{dt}{(1+t)^2}=1.
$$
Поэтому
 $$\lim\limits_{r\to \infty}\frac{\ln M(r,f)}{V_1(r)}=1.$$

Так как $V(r)\ge V_1(r)$, то 
$$
\varlimsup\limits_{r\to \infty}\frac{\ln M(r,f)}{V(r)}\le 1.
$$

Из того, что существует неограниченная последовательность $r_n$ такая, что $V_1(r_n)=V(r_n)$ следует, что
$$
\varlimsup\limits_{r\to \infty}\frac{\ln M(r,f)}{V(r)}=1.
$$
Теорема доказана.

\begin{center}

\end{center}
\textbf{А.Ф. Гришин (A. F. Grishin)}\\
Харьковский национальный университет
им. В. Н. Каразина,
\\ пл. Свободы, 4, 61022, Харьков, Украина\\
\textit{E-mail:} grishin@univer.kharkov.ua\\

\medskip
\noindent \textbf{Нгуен Ван Куинь (Nguyen Van Quynh)}\\
Харьковский национальный университет
им. В. Н. Каразина,
\\ пл. Свободы, 4, 61022, Харьков, Украина\\
\textit{E-mail:} quynh$\_$sonla032@yahoo.com; quynhsonla1988@gmail.com

\end{document}